\documentclass[reqno]{article}
\usepackage{amssymb, amsmath, amsthm, amsfonts, amscd}
\usepackage[english]{babel}

\usepackage{mathtools}
\usepackage{color}
\usepackage{hyperref}
\hypersetup{
    colorlinks=true,
    citecolor=red,
    filecolor=black,
    linkcolor=blue,
    urlcolor=black
}

\chardef\bslash=`\\ 

\newtheorem{thm}{Theorem}[section]

\newtheorem{lem}[thm]{Lemma}

\newtheorem{defn}{Definition}[section]

\newcommand{\N}{\mathbb{N}}

\newcommand{\ZZ}{\mathcal{Z}}
\newcommand{\Z}{\mathbb{Z}}
\newcommand{\Q}{\mathbb{Q}}
\newcommand{\R}{\mathbb{R}}
\newcommand{\C}{\mathbb{C}}
\newcommand{\Sp}{\mathbb{S}}
\newcommand{\T}{\mathbb{T}}

\newcommand{\ft}{\mathfrak{t}}

\def\a{\alpha }

\def\p{\pi}
\def\r{\rho}
\def\s{\sigma}

\def\D{\Delta}

\def\e{\varepsilon}

\def\.{\cdot }
\def\ra{\rightarrow}
\def\hra{\hookrightarrow}

\title{
\textsc{\textbf{Differentiable Rigidity for quasiperiodic cocycles in compact Lie groups}}\\
\author{Nikolaos Karaliolios \footnote{Affiliation: CNRS, IMJ-PRG,
UMR 7586, Univ. Paris Diderot, Sorbonne Paris Cit\'{e}, Sorbonne Universit\'{e}s, UPMC Univ Paris 06, F-75013, Paris, France.
} \footnote{email: nikolaos.karaliolios@imj-prg.fr} }
}

\begin{document}

\maketitle

\begin{abstract}
We study close-to-constants quasiperiodic cocycles in $\T ^{d} \times G$,
where $d \in \N ^{*} $ and $G$ is a compact Lie group, under the assumption that the rotation in
the basis satisfies a Diophantine condition. We prove differentiable rigidity for such cocycles: if such a cocycle
is measurably conjugate to a constant one satisfying a Diophantine condition with respect to the rotation,
then it is $C^{\infty}$-conjugate to it, and the K.A.M. scheme actually produces a conjugation. We also derive
a global differentiable rigidity theorem, assuming the convergence of the renormalization scheme for such
dynamical systems.
\end{abstract}  

\tableofcontents

\section{Introduction}
In the author's PhD thesis \cite{NKPhD}, the study of quasiperiodic cocycles in $\T ^{d} \times G$,
with $G$ a semisimple compact Lie group, over a Diophantine rotation and satisfying a closeness to constants
assumption was revisited. The basic reference for the subject is \cite{KrikAst}, where the corresponding
local density theorem is proved in the $C^{\infty}$ category by means of a K.A.M. scheme. The problem of
loss of periodicity (i.e. of conjugations of accumulating periods longer than $1$ in the presence of resonances) was settled outside the
iterative step of the scheme and made necessary the combination of the local almost quasi-reducibility theorem
with the reducibility theorem in a positive measure set of parameters in order to obtain a proof of
local almost reducibility. We were able to deal with this complication by proving a more efficient local
conjugation lemma, in which the phenomenon of longer periods is no longer present, and which can serve as an
iterative step of a K.A.M. scheme.

This improved scheme can be used in the proof of a local
differentiable rigidity theorem which improves the one obtained
in \cite{HY09}. Here,
$G$ is a semisimple compact Lie group.
\begin{thm} \label{Loc Dif Rig}
We suppose that $\a \in DC (\gamma , \tau ) $, and that $(\a , Ae^{F(\. )}) \in SW^{\infty} _{\a} (\T ^{d} ,G)$ is
close enough to $(\a ,A)$ so that the K.A.M. scheme can be initiated.
Here, $A\in G$ is a constant. Moreover, we suppose that there exists
$D(\. ) :\T ^{d} \ra G $, measurable, such that $Conj_{D(\. )}(\a , Ae^{F(\. )}) = (\a , A_{d}) $, where
$A_{d} \in DC _{\a} \subset G $. Then, the K.A.M. scheme can be made (with an appropriate adjustment of a
parameter) to produce only a finite number of resonances, and therefore to produce a $C^{\infty}$ conjugation.
In particular, $(\a , Ae^{F(\. )})$ is reducible.
\end{thm}
The Diophantine condition $DC$ is defined in def. \ref{def DC}, and condition
$DC_{\a }$ in \ref{Diophantine constants in G}.
The proof is carried out in $\T ^{d} \times SU(2)$, then
extrapolated to general compact Lie group $G$ by means of an appropriate
embedding $SU(2) \hra G$, so that no background on Lie group theory is
demanded for the greatest part of the note. Of course, $DC$ is
of full measure in $\T ^{d}$, and for every fixed $\a $  the condition
$DC_{\a }$ is of full Haar measure in $G$. The improvement in
comparison with the article cited above consists in the more
general algebraic context of our theorem, but mainly in the
fact that the smallness of the perturbation in our theorem
is not related with the measurable conjugation.

Subsequently, we briefly discuss the opposite phenomenon observed for "generic" cocycles smoothly conjugate to Liouville constant cocycles,
where the K.A.M. scheme produces an infinite number of resonances, and address the reader to \cite{NKInvDist}
for an improvement of the scheme which settles this problem.

Finally, we use the convergence of renormalization, as well as the measurable invariance of the degree
(see \cite{NKPhD}) in order to obtain a global differentiable rigidity theorem, without any assumtpion of
closeness to constants, but valid only for one-frequency cocycles ($d=1$).
\begin{thm} \label{Rigidity}
We suppose that $\a \in RDC $, and that $(\a , A(\. )) \in SW^{\infty} _{\a} (\T ,G)$ is
of degree $0$. Moreover, we suppose that there exists $D(\. ) :\T ^{d} \ra G $, measurable, such that
$Conj_{D(\. )}(\a , A(\. )) = (\a , A_{d}) $, where
$A_{d} \in DC _{\a} \subset G $. Then,  $(\a , A(\. ))$ is reducible.
\end{thm}
The assumption that $\deg (\a , A(\. )) =0 $, which we will not define here, assures that renormalization (see
\cite{Krik2001},\cite{AFK2010},\cite{AK2006},\cite{NKPhD}) converges to constants. This fact, combined
with the assumption $\a \in RDC$ (also of full measure in $\T$), implies that there exists
$\tilde{D}(\. ) \in C^{\infty}(\T ,G)$, such that $Conj_{\tilde{D}(\. )}(\a , A(\. ))$ satisfies the assumptions of
theorem \ref{Loc Dif Rig}. The condition $\deg (\a , A(\. )) =0$ is open and
dense in the $C^{\infty}$-topology, provided that $\a \in RDC$. We will not
come back to the proof of this theorem, since its details exceed the scope of
a short note.

\textbf{Acknowledgment}: The author would like to thank prof.
Rapha\"{e}l Krikorian, for bringing the question to his attention.
This work was partially supported by the ERC Starting Grant
"Quasiperiodic" and it appeared in a different form in
the author's PhD thesis.

\section{Facts from algebra and arithmetics}

\subsection{The group $SU(2)$}
The matrix group $G = SU(2)$ is the multiplicative group of unitary $2 \times 2$ matrices of determinant $1$.

Let us denote the matrix $S\in G$, $S=\begin{pmatrix}
z & w \\ 
-\bar{z} & \bar{w}
\end{pmatrix}
$, where $(z,w)\in \C^{2}$ and $|z|^{2}+|w|^{2}=1$, by $\{z,w\}_{G}$, and the subscript will be ommited unless
necessary. The manifold $G=SU(2)$ is thus naturally
identified with $\Sp ^{3}\subset \C^{2}$.
When coordinates in $\C ^{2}$ are fixed, the circle $\Sp ^{1}$ is naturally embedded in $G$ as the group of
diagonal matrices, which is a maximal torus (i.e. a maximal abelian subgroup) of $G$. The center of $G$,
noted by $Z_{G}$ is equal to $\{ \pm Id \} $.

The Lie algebra $g=su(2)$ is naturally isomorphic to $\R^3 \approx \R \times \C$ equipped with its vector and scalar product. It will be denoted by $g$. The element $s=\begin{bmatrix}
it & u \\ 
-\bar{u} & -it
\end{bmatrix}
$ will be denoted by $ \{t,u\}_{g} \in \R \times \C$.
Mappings with values in $g$ will be denoted by
\begin{equation*}
U(\. ) = \{ U_{t}(\. ), U_{z}(\. ) \}_{g} =  U_{t}(\. )h + U_{z}(\. ) j
\end{equation*}
in these coordinates, where $U_{t}(\. )$ is a real-valued and $U_{z}(\. )$ is a complex-valued function. The vectors $h,j,ij$ form an orthonormal and positively oriented basis for $su(2)$.

The adjoint action of $h\in su(2)$ on itself is pushed-forward to twice the vector product:
\begin{equation*}
ad_{ \{ 1,0 \} }.\{ 0,1 \} =  [ \{1,0\}, \{0,1\} ]  = 2\{ 0,i \}
= 2ij
\end{equation*}
plus cyclic permutations, and the Cartan-Killing form, normalized by
$\left\langle h,h^{\prime }\right\rangle =-\frac{1}{8\pi }tr(ad(h)\circ ad(h^{\prime }))$ is pushed-forward to the scalar
product of $\R^{3}$.
The preimages of the $Z_{G}$ in the maximal toral (i.e. abelian) algebra of
diagonal matrices are points of coordinates in the lattice $\pi \Z$.

The adjoint action of the group on its algebra is pushed-forward to the
action of $SO(3) \approx SU(2) / \pm Id $ on $\R \times \C$. In
particular, for diagonal matrices, of the form
$S = \exp (\{2\pi s,0\}_{g})$,
$Ad(S).\{ t,u\} =\{ t,e^{4i\pi s }u\} $.

\subsection{General compact groups}

The only fact that we will need from the theory of semisimple need groups is the decomposition of the Lie algebra,
$g'$, of such a group $G'$ in factors isomorphic to $g=su(2)$. If $g' $ is such an algebra, then for any
maximal abelian algebra $\ft $ (i.e. a subalgebra on which the restriction of the Lie bracket vanishes
identically), then the orthogonal complement of $\ft $ admits a decomposition
into a direct sum of mutually orthogonal $2$-dimensional real subspaces
$E_{\r }$, $\r \in \Delta _{+}$,\footnote{The finite set $\Delta _{+}$ is called the set of positive roots of $g'$,
and it is a subset of $\ft ^{*}$} for which the following holds. For each $E_{\r }$, there exists
$h_{\r } \in \ft$ such that $\R h_{\r } \oplus E_{\r } \approx su(2) $.
One can chose a vector $j _{\r }\in E_{\r }$, such that the vectors
$h_{\r} , j_{\r} , ij _{\r} $ form a basis of $\R h_{\r } \oplus E_{\r }$,
satisfying the same relations as the basis of $su(2)$.

There exists, moreover, $\tilde{\Delta} \subset \Delta _{+}$,
such that $(h_{\r })_{\tilde{\Delta}}$ is a basis for $\ft $, for
which there exists $P \in N ^{*} $ and $m_{\r , \r '} \in \N $
such that, for all $\r ' \in \Delta _{+} $,
\begin{equation} \label{lin sys res}
h_{\r '} = \frac{1}{P} \sum _{\r \in \Delta } m_{\r , \r '} h_{\r , \r '}
\end{equation}

The set of roots $\D _{+}$ satisfies the property that for all $\r \in \D _{+}$,
and for all $a \in \ft $, $[a,j_{\r}] = 2i\pi \r (a) j_{\r } = 
2i\pi a_{\r} j_{\r }$. The
adjoint action of $G$ on $\ft $ (and therefore on the vectors
$h_{\r }$) induces an action of $G$ on $\D _{+}$.

The Lie algebra $g$ is thus decomposed into
\begin{equation*}
\mathfrak{c} \oplus \ft \oplus \left( \oplus _{\r \in \D _{+}} E_{\r } \right) =
\mathfrak{c} \oplus \oplus _{\r \in \tilde{\D }} \R h_{\r }
\oplus  _{\r \in \D _{+}} \C j_{\r }
\end{equation*}
Here, $\mathfrak{c} $ is the center of the algebra, which is trivial if
the algebra is semisimple. This decomposition is referred to as "root space
decomposition with respect to the toral algebra $\ft $".

The lattice of preimages of $Z_{G }$ in $\ft$ will be denoted
by $\ZZ $.

\subsection{Calculus and Functional Spaces}

\subsubsection{Functional spaces}

We will consider the space $C^{\infty }(\T ^{d} ,g)$ equipped with the
standard maximum norms%
\begin{equation*}
\left\Vert U\right\Vert _{s} = \max_{0\leq \sigma \leq s} \max_{\T ^{d} }\left\vert \partial
^{\sigma }U(\. )\right\vert
\end{equation*}
for $s\geq 0$, and the Sobolev norms
\begin{equation*}
\left\Vert U\right\Vert _{H^{s}}^{2}=\sum_{k\in \Z^{d}}(1+|k|^{2})^{s}|\hat{U}(k)|^{2}
\end{equation*}
where $\hat{U}(k)=\int U(\. )e^{-2i\pi kx}$ are the Fourier coefficients
of $U(\. )$. The fact that the injections $H^{s +d/2}(\T ,g) \hra C^{s}(\T ^{d} ,g) $
and $C^{s}(\T ^{d},g) \hra H^{s}(\T ^{d},g)$ for all $s \geq 0$ are continuous is classical.

We will also use the convexity or Hadamard-Kolmogorov inequalities (see \cite{KolmIneq49}) ($U\in
C^{\infty }(\T ,g)$):
\begin{equation*}
\left\Vert U(\. )\right\Vert _{\sigma }\leq C _{s, \sigma }\left\Vert
U\right\Vert _{0}^{1-\sigma /s}\left\Vert U\right\Vert _{s}^{\sigma /s}
\label{hadamard}
\end{equation*}
for $0\leq \sigma \leq s$, and the inequalities concerning the composition of functions (see \cite{KrikAst}):
\begin{equation*} \label{composition}
\left\Vert \phi \circ (f+u)-\phi \circ f\right\Vert _{s}\leq
C_{s}\left\Vert \phi \right\Vert _{s+1}(1+\left\Vert f\right\Vert
_{0})^{s}(1+\left\Vert f\right\Vert _{s})\left\Vert u\right\Vert _{s}
\end{equation*}

We will use the truncation operators for mappings $\T ^{d} \ra g$ defined by
\begin{eqnarray*}
T_{N}f(\. ) &=&\sum_{|k|\leq N}\hat{f}(k)e^{2i\pi k\. } \\
\dot{T}_{N}f(\. ) &=&T_{N}f(\. )-\hat{f}(0) \\
R_{N}f(\. ) &=&\sum_{|k|>N}\hat{f}(k)e^{2i\pi k\. }
\end{eqnarray*}
These operators satisfy the estimates
\begin{eqnarray} \label{truncation est}
\left\Vert T_{N}f(\. )\right\Vert _{C^{s}} &\leq
&C_{s}N ^{d/2} \left\Vert f(\. )\right\Vert _{C^{s}} \\
\label{truncation rest}
\left\Vert R_{N}f(\. )\right\Vert _{C^{s}} &\leq &C_{s,s'} N^{s-s^{\prime }+d/2} \left\Vert
f(\. )\right\Vert _{C^{s^{\prime }}}
\end{eqnarray}
The Fourier spectrum of a function will be denoted by $\hat{\sigma}(f)=\{k\in \Z ^{d},~\hat{f}(k)\not=0\}$.

\subsection{Arithmetics, continued fraction expansion} \label{Arithmetics, continued fraction expansion}

The following notion is essential in K.A.M. theory. It is related with the
quantification of the closeness of rational numbers to certain classes of
irrational numbers.
\begin{defn} \label{def DC}
We will denote by $DC(\gamma ,\tau )$ the set of numbers $\a$ in $\T \setminus \Q$ such that for any $k\not=0$,
$|\a k| _{\Z }\geq \frac{\gamma ^{-1}}{|k|^{\tau }}$. Such numbers are called \textit{Diophantine}.
\end{defn}
The set $DC(\gamma ,\tau )$, for $\tau >2$ fixed and $\gamma \in \R _{+} ^{\ast}$ is of positive
Haar measure in $\T $. If we fix $\tau$ and let $\gamma$ run through the positive real numbers, we obtain
$\cup_{\gamma >0} DC(\gamma ,\tau )$ which is of full Haar measure. The numbers that do not satisfy any Diophantine
condition are called \textit{Liouvillean}. They form a residual set of $0$ Lebesgue measure.

This last following definition concerns the relation of the approximation of an irrational number with its continued
fractions representation.
\begin{defn} \label{def RDC}
We will denote by $RDC(\gamma ,\tau )$ is the set of  \textit{recurrent Diophantine} numbers, i.e. the $\a $ in
$\T \setminus \Q$ such that $G^{n}(\a)\in DC(\gamma ,\tau )$ for infinitely many $n$.
\end{defn}
Here, $G(\a)=\{\a^{-1}\}$ is the Gauss map ($\{ \. \} $ stands for
"fractional part"). The set $RDC$ is also of full measure, since the Gauss
map is ergodic with respect to a smooth measure.

In contexts where the parameters $\gamma $ and $\tau $ are not significant, they will be omitted in the notation of both sets.

\section{Cocycles in $\T ^{d} \times G$}

\subsection{The dynamics}

Let $\a \in \T ^{d} \equiv \R ^{d} / \Z ^{d} $, $d \in \N ^{*}$, be an irrational rotation, so that the
translation $x \mapsto x+\a \mod(\Z ^{d} ) $ is minimal and uniquely ergodic. The translation will sometimes
be denoted by $R_{\a } $.

If we also let $A(\. )\in C^{\infty}(\T ^{d} ,G)$, the couple $(\a ,A(\. ))$ acts on the fibered space
$\T ^{d} \times G \ra \T ^{d} $ defining a diffeomorphism by
\begin{equation*}
(\a ,A(\. )).(x,S)=(x+\a ,A(x).S), (x,S)\in \T ^{d} \times G
\end{equation*}
We will call such an action a \textit{quasiperiodic cocycle over }$R_{\a}$ (henceforth simply
a cocycle). The space of such actions is denoted by $SW_{\a }^{\infty}(\T ^{d},G)\subset
Diff^{\infty}(\T ^{d} \times G)$.
Most times we will abbreviate the notation to $SW_{\a }^{\infty}$. Cocycles are a class of fibered diffeomorphisms, since
fibers of $\T ^{d} \times G$ are mapped into fibers, and the mapping from one fiber to another in general depends on the base point.
The number $d \in \N ^{*}$ is the number of frequencies of the cocycle.

If we consider a representation of $G$ on a vector space $E $, the action of the cocycle can be also defined
on $\T ^{d} \times E $, simply by replacing $S$ by a vector in $E$ and multiplication in $G $ by the action.
The particular case which will be important in this article is the representation of $G $ on $g$,
and the resulting action of the cocycle on $\T ^{d} \times g$.

The $n$-th iterate of the action is given by
\begin{equation*}
(\a,A(\. ))^{n}.(x,S)=(n\a,A_{n}(\. )).(x,S)=(x+n\a 
,A_{n}(x).S)
\end{equation*}%
where $A_{n}(\. )$ represents the \textit{quasiperiodic product} of matrices equal to
\begin{equation*}
A_{n}(\. )= \begin{cases} A(\. +(n-1)\a)\cdots A(\. ) &, n >0 \\
Id &, n =0 \\
A^{\ast }(\. +n\a)\cdots A^{\ast }(\. -\a) &, n <0
\end{cases}
\end{equation*}

\subsection{Classes of cocycles with simple dynamics, conjugation}

The cocycle $(\a  ,A(\.  ))$ is called a constant cocycle if
$A(\. )=A\in G$ is a constant mapping. In that case,
the quasiperiodic product reduces to a simple product of matrices,
$(\a  ,A)^{n}=(n\a  ,A^{n})$.

The group $C^{\infty }(\T ^{d},G)\hra SW^{\infty }(\T ^{d},G)$ acts by \textit{dynamical conjugation}: Let
$B(\.  )\in C^{\infty }(\T ^{d},G)$ and $(\a  ,A(\. ))\in SW^{\infty }(\T ^{d},G)$. Then we define
\begin{eqnarray*}
Conj_{B(\.  )}.(\a  ,A(\.  )) &=&(0,B(\.  ))\circ (\a  ,A(\.  ))\circ (0,B(\.  ))^{-1} \\
&=& (\a  ,B(\.  +\a  ).A(\. ).B^{-1}(\.  ))
\end{eqnarray*}%
which is in fact a change of variables within each fiber of the product $\T ^{d} \times G$. The dynamics of
$Conj_{B(\. )}.(\a  ,A(\. ))$ and $(\a ,A(\. ))$ are essentially the same, since
\begin{equation*}
(Conj_{B(\.  )}.(\a  ,A(\.  )))^{n}=(n\a  ,B(\.  +n\a ).A_{n}(\.  ).B^{-1}(\. ))
\end{equation*}
\begin{defn}
Two cocycles $(\a ,A(\. ))$ and $(\a ,\tilde{A}(\.  ))$ in $SW^{\infty }_{\a }$ are \textit{conjugate} iff there
exists $B(\.  )\in C^{s}(\T ^{d},G)$ such that $(\a  ,\tilde{A}(\.  ))=Conj_{B(\.  )}.(\a  ,A(\.  ))$.
We will use the notation $(\a  ,A(\.  ))\sim (\a  ,\tilde{A}(\.  ))$
to state that the two cocycles are conjugate to each other.
\end{defn}
Since constant cocycles are a class for which dynamics can be analysed, we give the following definition.
\begin{defn} \label{defn a.r.}
A cocycle will be called \textit{reducible} iff it is conjugate to a constant.
\end{defn}
Due to the fact that not all cocycles are reducible (e.g. generic cocycles in $\T \times \Sp ^{1}$
over Liouvillean rotations, but also cocycles over
Diophantine rotations, even though this result is hard to obtain, see \cite{El2002a}, \cite{Krik2001})
we also need the following concept, which has proved to be crucial in the study of such dynamical systems.
\begin{defn} \label{def almost reducibility}
A cocycle $(\a ,A(\. ))$ is said to be \textit{almost reducible} if there exists a sequence of conjugations
$B_{n}(\. ) \in C^{\infty}$, such that $Conj_{B_{n}(\. )}.(\a ,A(\. ))$ becomes arbitrarily close to constants
in the $C^{\infty }$ topology, i.e. iff there exists $(A_{n})$, a sequence in $G$, such that
\begin{equation*}
A_{n}^{\ast } \left(  B_{n}(\. +\a )A(\. )B_{n}^{\ast }(\. ) \right) = e^{F_{n} (\. )} \overset{C^{\infty } }{\ra } Id
\end{equation*}
\end{defn}
This property will herein be established in a K.A.M. constructive way, making it possible to measure the rate of
convergence versus the explosion of the conjugations. Almost reducibility then comes along with obtaining that
\begin{equation*}
Ad (B_{n} (\. ) ) .F_{n} (\. ) = B_{n} (\. ) .F_{n} (\. ) . B_{n}^{*} (\. ) \overset{C^{\infty } }{\ra } 0
\end{equation*}
If this additional condition is satisfied, almost reducibility in the sense of
the definition above and almost reducibility in the sense that "the cocycle can be conjugated arbitrarily close
to reducible cocycles" are equivalent.

We can now recall the local almost reducibility theorem. It is used the proof of the local density theorem,
already proved in \cite{KrikAst}. We will follow the proof in \cite{NKPhD}, since it implies a useful corollary
(cor. \ref{preliminary reduction}) in a more direct manner than the previously existing proofs.
\begin{thm} \label{thm local a.r.}
Let $\a \in DC(\gamma , \tau ) \subset \T ^{d}, d\geq 1 $ and $G$ a semisimple compact Lie group. Then, there
exists $s_{0} \in \N ^{*} $ and $\epsilon >0$, such that if
$(\a , Ae^{F(\. )} )\in SW^{\infty}_{\a} (\T ^{d}, ) $ with $\| F(\. ) \| _{0} < \epsilon $
and $\| F(\. ) \| _{s_{0}} <1 $, $(\a , Ae^{F(\. )} )$ is almost reducible.
\end{thm}

\section{Almost reducibility} \label{sec a.r. and due}

In this section we present the basic points of the proof of the thm \ref{thm local a.r.}.
For the next paragraph, $G = SU(2)$, and the proof of the local conjugation lemma when $G$ is an arbitrary compact Lie group will be hinted in the next one.

\subsection{Local conjugation in $\T ^{d} \times SU(2)$}
Let $(\a , Ae^{F(\. )} ) =(\a , A_{1}e^{F_{1}(\. )} ) \in SW^{\infty} (\T ,G) $ be a cocycle over a
Diophantine rotation satisfying some smallness conditions to be made more precise later on. Without any loss of
generality, we can also suppose that $A = \{ e^{2i\pi a} , 0 \} $ is diagonal. The goal is to conjugate the cocycle
ever closer to constant cocyles by means of an iterative scheme. This
is obtained by iterating the following lemma,
for the detailed proof of which we refer to \cite{KrikAst},
\cite{El2002a} or \cite{NKPhD}. The following lemma is the
cornerstone of the procedure, since it represents one step of the
scheme. The rest of this paragraph is devoted to a summary of its
proof, for the sake of completeness.
\begin{lem} \label{loc conj lem}
Let $\a  \in DC(\gamma ,\tau )$ and $K\geq C \gamma N^{\tau }$. Let, also,
$(\a  ,Ae^{F(\. )})\in SW^{\infty }(\T _{d} ,G)$ with
\begin{equation*}
c_{1,0}KN ^{s_{0}} \e  _{0}<1
\end{equation*}%
for some $s_{0}\in \N ^{*}$ depending on $d, \gamma , \tau$, and where $\e  _{s}=\left\Vert F\right\Vert _{s}$.
Then, there exists a conjugation $G(\.  )\in C^{\infty }(\T ,G)$ such that
\begin{equation} \label{Conj eq}
G(\.  +\a  ).A.e^{F(\.  )}.G^{*}(\.  )=A^{\prime} e^{F' (\.  )}
\end{equation}
The mappings $G(\. )$ and $F'(\. )$ satisfy the following estimates
\begin{eqnarray*}
\left\Vert G(\.  )\right\Vert _{s} &\leq &  c_{1,s} (N^{s} + KN^{s+\tau +1/2}\e  _{0}) \\
\e  _{s}^{\prime } &\leq &c_{2,s}K^{2}N^{2\tau + d}(N^{s}\e  _{0}+\e _{s})\e  _{0}+
C_{s,s^{\prime }}N^{s-s^{\prime }+2\tau + d}\e _{s^{\prime }}
\end{eqnarray*}
\end{lem}

If we suppose that $Y(\. ) : \T \ra g $ can conjugate $(\a , Ae^{F(\. )} )$ to
$(\a , A'e^{F'(\. )} )$, with $\| F'(\. ) \| \ll \| F(\. ) \| $, then it must satisfy the functional
equation
\begin{equation*}
A^{*}e^{Y(\. +\a )}Ae^{F(\. )}e^{-Y(\. )}  =A^{*} A'e^{F'(\. )} 
\end{equation*}
Linearization of this equation under the assumption that all $C^{0}$ norms are smaller than $1$ gives
\begin{equation*} 
Ad(A^{*})Y(\. +\a ) +F(\. )-Y(\. )  = \exp ^{-1} (A^{*} A')
\end{equation*}
which we will write in coordinates, separating the diagonal from the non-diagonal part.

The equation for the diagonal coordinate reads $Y_{t } (\. +\a )-Y_{t } (\.)=-F_{t }(\. )$.
For reasons well known in K.A.M. theory, we have to truncate at an order $N $ to be determined by the parameters
of the problem and obtain a solution to the equation
\begin{equation*}
Y_{t } (\. +\a )-Y_{t } (\.)=- \dot{T}_{N} F_{t }(\. )
\end{equation*}
satisfying the estimate $\| Y_{t } (\.) \| _{s} \leq \gamma C_{s} N^{s + \tau + d/2} \| F_{t } (\.) \| _{0} $.
The rest satisfies the estimate of eq. \ref{truncation rest}. The mean value $\hat{F}_{t } (0) $ is an
\textit{obstruction} and will be integrated in $\exp ^{-1} (A^{*} A') $.

As for the equation concerning the non-diagonal part, it reads
\begin{equation} \label{untruncated local eq}
e^{-4i\pi a}Y_{z }(\. +\a )-Y_{z }(\. )=-F_{z}(\. )
\end{equation}
or, in the frequency domain,
\begin{equation} \label{local eq in Fourier}
(e^{2i\pi ( k\a  - 2a )}-1)\hat{Y}_{z  }(k)=-\hat{F}_{z }(k),~ k \in \Z ^{d}
\end{equation}
Therefore, the Fourier coefficient $\hat{F}_{z }(k_{r })$ cannot be eliminated with good estimates if
\begin{equation*}
| k_{r } \a  -2a | _{\Z } < K ^{-1}
\end{equation*}
for some $K>0 $ big enough. If $K = N^{\nu }$, with $\nu > \tau $,
then we know by \cite{El2002a} that, if such a $k_{r }$ exists
(called a \textit{resonant mode}) and satisfies $0\leq |k_{r} | \leq N $,
it is unique in $\{ k \in \Z ^{d} , |k- k_{r } | \leq 2N \}$.
Therefore, if we call $T^{k_{r }} _{2N} $ the truncation operator
projecting on the frequencies $ 0 < |k- k_{r } | \leq 2N$ if $k_{r }$
exists (or on $ |k| \leq N$ if it does not, but this easier and
follows from this one), the equation
\begin{equation*}
e^{-4i\pi a}Y_{z }(\. +\a )-Y_{z }(\. )=-T^{k_{r }} _{2N} F_{z}(\. )
\end{equation*}
can be solved and the solution satisfies
$\| Y_{z } (\.) \| _{s} \leq C_{s} N^{s + \nu + d/2} \| F_{z } (\.) \| _{0} $.
We will define the rest operator by projection on modes satisfying
$|k- k_{r } | > 2N$.

In total, the equation that can be solved with good estimates is
\begin{equation*}
Ad(A^{*})Y(\. +\a ) -Y(\. )  = - F(\. ) +
\{ \hat{F}_{t } (0), \hat{F} _{z} (k_{r }) e^{2i\pi k_{r } \.} \} +
\{ R_{N} F_{t} (\. ) , R^{k_{r }} _{2N} F_{z} (\. ) \}
\end{equation*}
with $\| Y(\. ) \|_{s} \leq C_{s} N^{s + \nu + 1/2 } \| F(\. ) \| _{0} $. Under the smallness assumptions of the
hypothesis, the linearization error is small and the conjugation thus constructed satisfies
\begin{equation*}
e^{Y(\. +\a )}Ae^{F(\. )}e^{-Y(\. )}  =
\{ e^{2i\pi (a + \hat{F}_{t } (0) )} ,0 \} _{G} . e^{\{ 0, \hat{F} _{z} (k_{r }) e^{2i\pi k_{r } \.} \} _{g} }
e^{\tilde{F} (\. )} 
\end{equation*}
with $\tilde{F} (\. ) $ a "quadratic" term. We remark that, a priori, the obstruction
$\{0, \hat{F} _{z} (k_{r }) e^{2i\pi k_{r } \.} \}$ is of the order of the initial perturbation and therefore
what we called $\exp ^{-1} (A^{*} A') $ is not constant in the presence of resonant modes.

If $k_{r }$ exists and is non-zero, then application the lemma cannot
be iterated. On the other hand, the conjugation
$B(\. ) = \{ e^{- 2i\pi k_{r } \. /2} ,0 \} $ is such that, if we call
$F '' (\. ) = Ad(B(\. )).F(\. ) = \{ F_{t}(\. ) , e^{-2i\pi _{k_{r} \.} } F_{z} (\. ) \} $, and $Y'(\. ) = Ad(B(\. )). Y(\. )$,
and $A'' = B(\a ) A = \{ e^{2i\pi (a - k_{r} \a /2)} ,0 \}$,
they satisfy the equation
\begin{eqnarray*}
Ad((A'')^{*})Y'(\. +\a ) -Y'(\. ) +F_{1}^{\prime \prime}(\. ) &=&
\{ \hat{F}_{t } (0), \hat{F} _{z} (k_{r })  \} +
\{ R_{2N} F_{t} (\. ) , e^{- 2i\pi k_{r } \.} R^{k_{r }} _{2N} F_{z} (\. ) \} \\
&=& \{ \hat{F}_{t }^{\prime \prime} (0), \hat{F} _{z}^{\prime \prime} (0)  \} +
\{ R_{2N} F_{t}^{\prime \prime} (\. ) , \tilde{R}^{k_{r }} _{2N} F_{z} (\. ) \}
\end{eqnarray*}
where $ \tilde{R}^{k_{r }} _{2N}$ is a dis-centred rest operator, whose spectral support is outside $[-N,N]^{d}\cap \Z ^{d} $,
and can therefore be estimated like a classical rest operator
$R_{N}$. The equation for primed variables can be
obtained from eq. \ref{untruncated local eq} by applying $Ad(B(\. ))$ and using that $B(\. )$ is a
morphism and commutes with $A$. The passage from one equation to the other is equivalent to the fact that
\begin{eqnarray*}
Conj_{B(\. )} (\a , A.\exp (\{ \hat{F}_{t } (0), \hat{F} _{z} (k_{r }) e^{ 2i\pi k_{r } \.}) \} ) &=&
(\a , B(\a ) A.\exp (\{ \hat{F}_{t } (0), \hat{F} _{z} (k_{r })  \} ) \\
&=& (\a , \tilde{A})
\end{eqnarray*}
that is, $B(\. )$ reduces the initial constant perturbed by the obstructions to a cocycle close to $(\a , \pm Id )$.
There is a slight complication, as $B(\. )$ may be $2$-periodic (if
$\{ i\pi k_{r} , 0 \} \in g$ is a preimage of $-Id$). If it is so,
we can conjugate a second time with a minimal geodesic
$C(\. ): 2 \T \ra G$ such that $C(1) = -Id \in Z_{G} $ and commuting
with $\tilde{A}$. The cocycle that we obtain in this way is
$1$-periodic and close to $(\a , \{ e^{i\p \a} ,0 \} _{G} )$, and
the conjugation is also $1$-periodic.

Summing up, if we call $G(\. ) = C(\. ) B(\. )e^{Y(\. )}$ and
$A' = C(\a ) \tilde{A}$, there exists $F'(\. )$ satisfying the
estimates of the lemma and such that eq. \ref{Conj eq} is verified.

\subsection{Local conjugation in general Lie groups}

Local conjugation in a compact group $G$ is in fact a vector-space
version of the lemma of the previous paragraph. The statement
is the same, just replacing $SU(2)$ by a compact group $G$, and
the steps of the proof are the same. Firstly, solution of the linear
equation, which is done in $\mathfrak{c} \oplus \ft $ as in the
diagonal coordinates here above, and in each $E_{\r}$ as in
eq. \ref{local eq in Fourier}. This procedure produces a "vector"
\footnote{Some entries may be $\emptyset$, if there is no resonant
mode in the corresponding eigenspace, or the corresponding frequency
in $\Z $}
$(k_{\r })_{\r \in \D _{+}} \in \bigsqcup _{\D _{+}} \Z $ of resonant
modes, which are reduced in the second step, where we construct a
vector $H$ such that
\begin{eqnarray*}
Conj_{\exp (H \. )} (\a , A.\exp (\{ Ob F(\. ) \} ) &=&
(\a , e^{H\a } A.\exp (\{ \hat{F}_{\ft \oplus \mathfrak{c} } (0)+ 
Ad (e^{H\. }). Ob F(\. )  \} ) \\
&=& (\a , \tilde{A})
\end{eqnarray*}
Here, $Ob$ stands for projection on the resonant modes,
$\sum _{\r \in \D _{+}} \hat{F} (k_{\r })e^{2i\pi \.} j_{\r } $.
The vector $J \in \frac{1}{P} \ZZ $
is constructed by solving the linear system of eq.
\ref{lin sys res} for the vector $(k_{\r })_{\r \in \D _{+}} $.
Finally, since the mapping $\exp (J \. )$ may be
$\frac{1}{P \# Z_{G} }$-periodic, we post-conjugate with $C(\. )$,
a minimal geodesic connecting the $Id$ with $\exp J$ (which measures
the failure of $\exp (J \. )$ to be $1$-periodic) and commuting with
$\tilde{A}$. We thus obtain $A'$ which is close to an element of
$\a \frac{1}{P} \ZZ $.\footnote{Intersected with the ball in $g$
where $\exp$ is bijective.}
\subsection{The K.A.M. scheme} Lemma \ref{loc conj lem} can serve as the step of a K.A.M. scheme, with the following
standard choice of parameters: $N_{n+1} = N_{n}^{1+\s } = N^{(1+\s )^{n-1}}$, where $N=N_{1}$ is big enough and
$0<\s <1$, and $K_{n} = N_{n}^{\nu }$, for some $\nu > \tau$. If we suppose that $(\a ,A_{n}e^{F_{n}(\. )})$
satisfies the hypotheses of lemma \ref{loc conj lem} for the corresponding parameters, then we obtain
a mapping $G_{n}(\. ) = C_{n}(\. ).B_{n}(\. ) e^{Y_{n}(\. )} $ that conjugates it to $(\a ,A_{n+1}e^{F_{n+1}(\. )})$,
and we use the notation $\e _{n,s} = \| F_{n} \|_{s}$.

If we suppose that the initial perturbation is small in small norm: $\e_{1,0} < \epsilon <1$, and not big in
some bigger norm: $\e _{1, s_{0}}<1$, where $\epsilon$ and $s_{0}$ depend on the choice of parameters,
then we can prove (see \cite{NKPhD} and, through that, \cite{FK2009}), that the scheme can be iterated,
and moreover
\begin{eqnarray*}
\e _{n,s} &=& O(N_{n}^{-\infty}) \text{ for every fixed } s \text{ and }\\
\| G_{n} \| _{s} &=& O (N_{n}^{s+\lambda}) \text{ for every } s \text{ and some fixed } \lambda >0
\end{eqnarray*}
We say that the norms of perturbations decay exponentially, while conjugations grow polynomially.

\subsection{A "K.A.M. normal form"}
The product of conjugations $H_{n} = G_{n} \cdots G_{1} $, which by construction satisfies
$Conj_{H_{n}} (\a ,A_{1} e^{F_{1}}) = (\a ,A_{n+1} e^{F_{n+1}})$, is not expected to converge. In fact,
it converges iff $B_{n} (\. ) \equiv Id $, except for a finite number of steps.
Anyhow, we can obtain a "\textit{K.A.M. normal form}" for cocycles close to constants
\begin{lem} \label{preliminary reduction}
Let the hypotheses of theorem \ref{thm local a.r.} hold. Then, there exists $D(\. ) \in C^{\infty } (\T ,G)$
such that, if we call $Conj_{D(\. )} (\a ,A e^{F(\. )}) = (\a ,A' e^{F'(\. )})$, then the K.A.M. scheme applied to
$(\a ,A' e^{F'(\. )})$ for the same choice of parameters consists only in the reduction of resonant modes: The
resulting conjugation $H_{n_{i}}(\. )$ has the form
$\prod \nolimits _{n_{i}} ^{1} C_{n_{i}} (\. ). B_{n_{i}} (\. ) $,
where $\{ n_{i} \}$ are the steps in which reduction of a resonant
mode took place.
\end{lem}
The proof, for which we refer the reader to \cite{NKInvDist}, is short and uses the fast convergence of the scheme.
For a cocycle in normal form, we relabel the indexes as $(\a ,A_{n_{i}}e^{F_{n_{i}}}) = (\a ,A_{i}e^{F_{i}})$.

\section{Proof of theorem \ref{Loc Dif Rig}}

In this section, we use the K.A.M. scheme outlined above. We begin by a brief study
of the rigidity of conjugation between constant cocycles in general compact
groups, which is to be compared with section \textit{2.5.e} of
\cite{KrikAst}, and shows that conjugation between constant cocycles is
rigid, with no assumptions on the arithmetics.

\subsection{The toy case}

Let $B:\mathbb{T}\rightarrow G$ be a measurable mapping, $\a  \in \T ^{d}$
a minimal translation, and $C_{1},C_{2}\in G$ such that%
\begin{equation*}
B(\.  +\a  )C_{1}B^{\ast }(\.  )=C_{2}
\end{equation*}%
By composing $B$ with a constant if necessary, we can suppose that $C_{1}$
and $C_{2}$ are on the same maximal torus. If, for simplicity in notation we
identify $G$ and $Inn(g) \approx G / Z_{G} $, the group which acts by the adjoint action on $g$,
even though the identification is not accurate, we find that%
\begin{equation*}
e^{2i\pi k\a  }\hat{B}(k)C_{1}=C_{2}\hat{B}(k)
\end{equation*}%
so that, if $e^{2i\pi c^{(i)}_{\r}}$ are the eigenvalues of the adjoint action of $C_{i}$, the equation%
\begin{equation*}
\langle e^{2i\pi (k\a  +c_{\r  }^{(1)})}\hat{B}(k)j_{\r  },j_{\r 
^{\prime }}\rangle =\langle \hat{B}(k)j_{\r  },e^{2i\pi c_{\r 
}^{(2)}}j_{\r  ^{\prime }}\rangle
\end{equation*}%
has at most one non-zero solution in $k$ for each pair of roots $\r  $ and $%
\r  ^{\prime }$ for which there exists $Z\in G$ such that $\r  ^{\prime
}=Z.\r  $. Using a similar argument for directions in the torus, we find
that $B:\mathbb{T}\rightarrow Inn(g)$, and therefore $B:\mathbb{T}%
\rightarrow G$, has a finite support in the frequency space, and therefore
it is $C^{\infty }$. Finally, it can be easily proved that $B(\. )$ can
be written as a product of two commuting morphisms, which do not take
values in the same maximal torus.

\subsection{Proof of theorem \ref{Loc Dif Rig}}

In order to simplify the proof and to avoid the phenomena
related to loss of periodicity, we consider a unitary representation of $G$
and we suppose firstly that $D(\. +\a )A_{1}e^{F_{1}(\. )}D^{\ast }(\. ) = A_{d}$
with $D(\.  ): \T \ra G\hra U(w)$ a measurable mapping and $F_{1}(\.  )$
small enough so that the reduction scheme can be applied. The K.A.M. scheme
of the previous section, applied in the simpler algebraic context of $U(w)$,
produces the sequence of conjugations
$H_{n_{i}}(\. ) = H_{i}(\. ) = B_{n_{i}}(\. ) \cdots B_{n_{1}}(\. )$.
This product converges if, and only if, it is finite, and satisfies
\begin{equation*}
D(\. +\a )H_{i}^{\ast }(\. +\a )A_{i}e^{F_{i}(\. )}H_{i}(\. )D^{\ast }(\. )=A_{d}
\end{equation*}%
or, by introducing some obvious notation%
\begin{equation*}
D_{i}(\.  +\a  )A_{i}e^{F_{i}(\.  )}D_{i}^{\ast }(\.  )
=A_{d}
\end{equation*}
Finally, by post-conjugating with a constant we can assume that
$A_{i}$ and $A_{d}$ are in the same maximal torus.

\subsubsection*{The case $G=SU(2)$}
Let us firstly examine the simpler case where $G=SU(2)$, and suppose
that $H_{i}(\. )$ diverges.If we call
$D_{i}^{j,j}(\. ) = \langle D_{i}(\. )j,j\rangle : \T^{d} \ra \C $, we have
\begin{equation*}
(e^{2i\pi (k\a  +a_{i}-a_{d})}-1)e^{2i\pi a_{d}}\hat{D}_{i}^{j,j}(k)
= O(\e  _{n,0})
\end{equation*}
where $O(\e  _{n,0})$ is bounded independently of $k$. The divergence of
$D_{i}(\.  )$ in $L^{2}$ implies that $|a_{i}|_{\Z }\leq K_{n_{i}}^{-1}=N_{n_{i}}^{-\nu }$,
where $A_{i}= \{ \exp (2i \p a_{i}) , 0 \} _{SU(2)}$. Therefore,
\begin{equation*}
|k\a  +a_{i}-a_{d}|_{\Z }\geq K_{n_{i}} , \forall
0<|k|\leq (2\tilde{\gamma}^{-1})^{1/\tau }N_{n_{i}}^{^{\nu /\tau }}
=N_{i}^{\prime }
\end{equation*}
provided that%
\begin{equation} \label{dioph csts gr def}
|k\a  -a_{d}|_{\Z }\geq \dfrac{\tilde{\gamma} ^{-1}}{ |k| ^{\tau } }
\end{equation}
for $k \in \Z ^{\ast } $, i.e. iff $a_{d} \in DC_{\a } (\tilde{\gamma} ,
\tau )$. This motivates the following definition.
\begin{defn} \label{Diophantine constants in G}
A constant $A = \exp (a) \in G$ is diophantine with respect to $\a $ iff
for all roots $\r $ of a root-space decomposition
with respect to a torus passing by $A$ we have $\r (a) \in DC_{\a }
(\tilde{\gamma} ,\tilde{\tau } )$, i.e. if eq. \ref{dioph csts gr def} is
satisfied the constants $\gamma ,\tilde{\tau }$. By abuse of notation
we will write $A \in DC_{\a } (\tilde{\gamma} , \tilde{\tau } )$.
\end{defn}
Since $\nu >\tau $, we find that $N_{i}^{\prime }/N_{n_{i}}$ goes to
infinity, and therefore, for $i$ big enough the spectral support of $%
H_{i}(\.  )= D_{i}(\.  )D^{\ast }(\.  )$ (i.e. of the diverging sequence of
conjugations) is contained in
$[- N_{i}^{\prime } , N_{i}^{\prime }]^{d} \subset \Z ^{d} $.
In a similar way, and with some obvious notation, we find that
\begin{equation*}
(e^{2i\pi (k\a  +a_{d})}-1)\hat{D}_{i}^{h,j}(k) = O(\e  _{n,0})
\end{equation*}
Consequently,%
\begin{equation*}
\dot{T}_{N_{i}^{\prime }}D_{i}(\.  )=O_{L^{2}}(\e  _{n,0})
\end{equation*}%
On the other hand, since $\s (H_{i} (\. )) \subset [-2N_{n_{i}},2N_{n_{i}}]^{d} $ and since
$N_{n_{i}} \ll N_{i}^{\prime} $.
\begin{eqnarray*}
T_{N_{i}^{\prime }/2}D(\.  ) &=&T_{N_{i}^{\prime }/2}[D_{i}(\. )D_{i}^{\ast }(\.  )D(\.  )] \\
&=&T_{N_{i}^{\prime }/2}[T_{N_{i}^{\prime }}(D_{i}(\.  ))~H_{i}^{\ast }(\.  )] \\
&=&C_{i}H_{n}^{\ast }(\.  ) +O_{L^{2}}(\e  _{n,0})
\end{eqnarray*}

Since for $n$ big enough $\left\Vert D(\.  )-T_{N_{i}^{\prime }}D(\.  )\right\Vert _{L^{2}}$ is small and
$D_{i}^{\ast }(\.  )D(\.  )$ takes values in $Inn(su(2))\approx SO(3)$, we can assume that $C_{i}\in L(su(2))$
is bounded away from $0$, say
\begin{equation*}
|C_{i}|>\frac{1}{2}
\end{equation*}%
in operator norm. Since $H_{i}^{\ast }(\.  )$ diverges in $L^{2}$, we
reach a contradiction.
\subsubsection*{The case of a general compact Lie group}
The case of a general compact group is hardly more complicated. If the
sequence of conjugations diverges, there exists a root $\r  $ such that $%
|a_{\r  }^{(i)}|_{\Z }\leq K_{n_{i}}^{-1}$, where $a_{\r} = \r (a)$.
If we fix such a root $\r  $, we find that for any other positive root
$\r  ^{\prime }$,
\begin{eqnarray*}
(e^{2i\pi (k\a  +a_{\r  }^{(n)}-a_{\r  }^{(d)})}-1)e^{2i\pi
a_{\r  }^{(d)}}\hat{D}_{i}^{j_{\r  },j_{\r  ^{\prime }}}(k)
&=&O(\e  _{n,0}) \\
(e^{2i\pi (k\a  +a_{d})}-1)\hat{D}_{i}^{h_{\r  ^{\prime
}},j_{\r  }} (k)&=&O(\e  _{n,0})
\end{eqnarray*}%
and, in the same way as before,
\begin{equation*}
\dot{T}_{N_{i}^{\prime }}D_{i}(\.  ).j_{\r  }=O_{L^{2}}(\e  _{n_{i},0})
\end{equation*}
and since%
\begin{eqnarray*}
T_{N_{i}^{\prime }}D(\.  ).j_{\r  } &=&T_{N_{i}^{\prime }}[D_{i}(\. 
)D_{i}^{\ast }(\.  )D(\.  ).j_{\r  }] \\
&=&T_{N_{i}^{\prime }}[T_{2N_{i}^{\prime }}(D_{i}(\.  ))~D_{i}^{\ast
}(\.  )B(\.  ).j_{\r  }] \\
&=&C_{i}D_{i}^{\ast }(\.  )D(\.  ).j_{\r  }+O_{L^{2}}(\e  _{n_{i},0})
\end{eqnarray*}%
Since $D(\.  )$ is an isometry, we find that $C_{i}\in L(g)$ is bounded away from $0$, say
\begin{equation*}
|C_{i}|>\frac{1}{2}
\end{equation*}%
in operator norm. Now, $D_{i}^{\ast }(\.  )D(\.  ).j_{\r  }$ diverges in $L^{2}$ as $j_{\r  }$ does not commute
with the reduction of resonant modes. This is due to the fact that,  by construction of $h_{i}$ and by
the choice of the root $\r $
\begin{equation*}
[ h_{i} , j_{\r }] =2i \pi \dfrac{k^{\prime}_{\r }}{D}
\end{equation*}
infinitely often, with $k^{\prime}_{\r } \ra \infty$. Thus, the hypothesis that the product
of conjugation diverges leads us to a contradiction.

Finally, we observe that if $U(\.  )$ is small enough so that the K.A.M.
scheme can be applied, the diophantine condition on $A_{d}$ becomes
irrelevant. If we suppose that $A_{d}\in DC_{\a  }(\tilde{\gamma},\tau
^{\prime })$ with $\tau ^{\prime }>\tau $, then, after a finite number of
iterations of the scheme, $F_{i}(\.  )$ is small enough so that the scheme
can be initiated if we place $\a  $ in $DC(\gamma ,\tau ^{\prime })$, and
the argument presented above remains valid, and this concludes the proof of the theorem in its full generality.

\subsection{Reducibility to a Liouvillean constant}

A corollary of this proof is, in fact, the optimality of the scheme in the
orbits of Diophantine constant cocycles.\ By its construction, the scheme
converges in the smooth category if, and only if, it converges in $L^{2}$,
and the proof implies that if a measurable conjugation to such a constant
exists, then the scheme converges toward it, eventually modulo a conjugation between constant cocycles.

\bigskip

On the other hand, the transposed argument shows that the K.A.M. scheme is highly
non-optimal if the dynamics in the fibres are Liouvillean. More precisely,
we let $(\a  ,Ae^{F(\.  )})\in SW^{\infty }(\mathbb{T},SU(2))$ be
smoothly conjugate to $(\a  ,A_{L})$, where $A_{L}$ is a Liouvillean
constant in $SU(2)$. Application of the scheme produces a sequence of
conjugations $F_{n}(\.  )$ and a sequence of cocycles $(\a 
,A_{n}e^{F_{n}(\.  )})$ such that%
\begin{equation*}
Ae^{F(\.  )}=H_{n}(\.  +\a  )A_{n}e^{F_{n}(\.  )}H_{n}^{\ast
}(\.  )
\end{equation*}%
where $F_{n}(\.  )\rightarrow 0$ exponentially fast. If we suppose that
the sequence of conjugations converges, we find that in the limit%
\begin{equation*}
A_{L}=\tilde{H}(\.  +\a  )A_{\infty }\tilde{H}^{\ast }(\.  )
\end{equation*}%
where $A_{\infty }$ (which we suppose diagonal, just as $A_{L}$) is the
limit of $A_{n}$. Since $A_{L}$ is non-resonant, $\tilde{H}(\.  )$ is a
torus morphism, so that $A_{\infty }$ is Liouvillean itself.
Since, now, $F_{n}(\.  )\rightarrow 0$ exponentially fast and, for $n$ big enough,
$A_{n+1}=A_{n}\exp (\hat{F}_{n}(0))$, we can rewrite $A_{n}e^{F_{n}(\.  )}$ as $A_{\infty }e^{\tilde{%
F}_{n}(\.  )}$ where still $\tilde{F}_{n}(\.  )\rightarrow 0$
exponentially fast. Since $A_{\infty }$ is Liouvillean, for any $l\in 
\mathbb{N}$, there exists $k_{l}$ such that%
\begin{equation*}
|a_{\infty }-k_{l}\a  |<\frac{1}{|k_{l}|^{l}}
\end{equation*}%
Therefore the Fourier mode $k_{l}$ is a resonance for the scheme at the $n$%
-th step provided that%
\begin{eqnarray*}
\frac{1}{|k_{l}|^{l}} &<&N_{n}^{-\nu  } \\
|k_{l}| &<&N_{n}
\end{eqnarray*}%
or equivalently%
\begin{equation*}
N_{n}^{\nu /l}<|k_{l}|<N_{n}
\end{equation*}%
Therefore, since for $l > \nu $ a reduction of resonant mode must take place,
which contradicts the hypothesis that $H_{n}(\.  )$ converges. Since $l$ can
be chosen arbitrarily big, no choice of $\nu $ can make the scheme converge.

Strictly speaking, this phenomenon appears for "generic" reducible
cocycles. Genericity here is in the product topology
$G \times C^{\infty} (\T ^{d},G)$, which is very far from the one
induced by $SW^{\infty}$. Genericity in this sense comes from the
demand that resonant modes be non-$0$ infinitely often, in order
to avoid trivialities. A converging K.A.M. scheme can nonetheless
be defined, but we refer the reader to \cite{NKInvDist} for its
construction.

\bibliography{aomsample}

\providecommand{\bysame}{\leavevmode\hbox to3em{\hrulefill}\thinspace}
\providecommand{\noopsort}[1]{}
\providecommand{\mr}[1]{\href{http://www.ams.org/mathscinet-getitem?mr=#1}{MR~#1}}
\providecommand{\zbl}[1]{\href{http://www.zentralblatt-math.org/zmath/en/search/?q=an:#1}{Zbl~#1}}
\providecommand{\jfm}[1]{\href{http://www.emis.de/cgi-bin/JFM-item?#1}{JFM~#1}}
\providecommand{\arxiv}[1]{\href{http://www.arxiv.org/abs/#1}{arXiv~#1}}
\providecommand{\doi}[1]{\url{http://dx.doi.org/#1}}
\providecommand{\MR}{\relax\ifhmode\unskip\space\fi MR }
\providecommand{\MRhref}[2]{%
  \href{http://www.ams.org/mathscinet-getitem?mr=#1}{#2}
}
\providecommand{\href}[2]{#2}
\begin{thebibliography}{AFK01}

\bibitem[AFK01]{AFK2010}
\bgroup\scshape{}A.~Arthur\egroup{}, \bgroup\scshape{}B.~Fayad\egroup{}, and
  \bgroup\scshape{}R.~Krikorian\egroup{}, A {KAM} scheme for ${SL}(2,{R})$
  cocycles with {L}iouvillean frequences,  \emph{arXiv} (2001).

\bibitem[AK06]{AK2006}
\bgroup\scshape{}A.~Avila\egroup{} and \bgroup\scshape{}R.~Krikorian\egroup{},
  Reducibility or nonuniform hyperbolicity for quasiperiodic {S}chr\"odinger
  cocycles,  \emph{Ann. of Math. (2)} \textbf{164} (2006), 911--940.
  \mr{2259248 (2008h:81044)}.

\bibitem[Eli02]{El2002a}
\bgroup\scshape{}L.~H. Eliasson\egroup{}, Ergodic skew-systems on ${T}^d \times
  {SO}(3,{R})$,  \emph{Ergodic Theory and Dynamical Systems} \textbf{22}
  (2002), 1429--1449.

\bibitem[FK09]{FK2009}
\bgroup\scshape{}B.~Fayad\egroup{} and \bgroup\scshape{}R.~Krikorian\egroup{},
  Rigidity results for quasiperiodic {${\rm SL}(2,\Bbb R)$}-cocycles,  \emph{J.
  Mod. Dyn.} \textbf{3} (2009), 497--510. \mr{2587083 (2011f:37046)}.

\bibitem[HY09]{HY09}
\bgroup\scshape{}X.~Hou\egroup{} and \bgroup\scshape{}J.~You\egroup{}, Local
  rigidity of reducibility of analytic quasi-periodic cocycles on {${\rm
  U}(n)$},  \emph{Discrete Contin. Dyn. Syst.} \textbf{24} (2009), 441--454.
  \mr{2486584 (2010e:37024)}.  \doi{10.3934/dcds.2009.24.441}.

\bibitem[Kar13]{NKPhD}
\bgroup\scshape{}N.~Karaliolios\egroup{}, \emph{Aspects globaux de la
  r\'{e}ductibilit\'{e} des cocycles quasi-p\'{e}riodiques \`{a} valeurs dans
  des groupes de Lie compacts semi-simples}, 2013. Available at
  \url{http://tel.archives-ouvertes.fr/docs/00/77/79/11/PDF/these.pdf}.

\bibitem[Kar14]{NKInvDist}
\bgroup\scshape{}N.~Karaliolios\egroup{}, \emph{Invariant distributions for
  quasiperiodic cocycles in $\T ^{d} \times SU(2)$, in preparation}, 2014.
  \arxiv{1407.4763}.

\bibitem[Kol49]{KolmIneq49}
\bgroup\scshape{}A.~Kolmogoroff\egroup{}, On inequalities between the upper
  bounds of the successive derivatives of an arbitrary function on an infinite
  interval,  \emph{Amer. Math. Soc. Translation} \textbf{1949} (1949), 19.
  \mr{0031009 (11,86d)}.

\bibitem[Kri99]{KrikAst}
\bgroup\scshape{}R.~Krikorian\egroup{}, \emph{R\'eductibilit\'e des syst\`emes
  produits-crois\'es \`a valeurs dans des groupes compacts}, no. 259, 1999.
  \mr{1732061 (2001f:37030)}.

\bibitem[Kri01]{Krik2001}
\bgroup\scshape{}R.~Krikorian\egroup{}, Global density of reducible
  quasi-periodic cocycles on ${T}^1 \times {SU}(2)$,  \emph{Annals of
  Mathematics} \textbf{154} (2001), 269--326.

\end{thebibliography}
\bibliographystyle{aomalpha}

\end{document}